\documentstyle[12pt,epsf]{article}

\def\text{\hbox} 
\def\RR{{ I\!\!R}}

\def\wL0{{\hbox{\tilde L}_0}}

\newcommand{\cul}[2]{\hbox{\raisebox{-1.5ex}{$\stackrel{\!\!#1}{
\stackrel{\:\mathop{\kern 0pt\hbox to 10mm{\rightarrowfill}}\:}
{\scriptstyle \!\!#2}}$}}}

\newcommand{\stackunder}[2]{#2\limits_{#1}}

\newtheorem{lemma}{Lemma}[section]

\begin{document}

\begin{center} {\large\bf Domain of attraction of 
the quasi-stationary distributions for the Ornstein-Uhlenbeck process.}

{\scriptsize by} 
{\scriptsize\bf MANUEL LLADSER} 
\footnote{\label{tb} Departamento de Ingenier\'{\i}a Matem\'{a}tica, 
Universidad de Chile,
Casilla 170, Correo 3, Santiago 3, Chile.}
\footnote{e-mail:  lladser@dim.uchile.cl} and 
{\scriptsize\bf JAIME SAN MARTIN}
$^1$ \footnote{e-mail:  jsanmart@dim.uchile.cl}

\end{center}

\begin{abstract}
Let $(X_t)$ be a one-dimensional Ornstein-Uhlenbeck process with initial 
density function $f:\RR_+\to\RR_+$, which is a regularly varying function
with exponent $-(1+\eta), \ \eta \in(0,1)$.
We prove the existence of a 
probability measure $\nu$ with a Lebesgue density, 
depending on $\eta$, such that for every $A\in{\cal{B}}(\RR_+)$:
\[\lim_{t\to\infty}P_f(X_t\in A \mid T_0^X>t)=\nu(A).\]
\end{abstract}


\section{Introduction}

Let $\Omega= {\cal C}([0,\infty),\RR)$ be the space of real continuous
functions, and ${\cal F}$ the standard Borel $\sigma$-field on
$\Omega$. For a probability measure $\mu$ on $(\RR,{\cal{B}}(\RR))$
we denote by $P_\mu$ the probability measure on $(\Omega,{\cal F})$
such that $B_t(w)=w(t)$ is a Brownian Motion with initial
distribution $\mu$. If $\mu=\delta_x$ is the Dirac mass at $x \in \RR$
we denote $P_x$ instead of $P_{\delta_x}$. Similarly, if $\mu$
has a density $f$ we use the notation $P_f$ instead
of $P_\mu$.

Consider a one-dimensional diffusion process $(X_t)$ which in 
differential form may be written as
\begin{equation} \label{gede} 
dX_t=dB_t-\alpha (X_t)dt, \; X_0=B_0,  
\end{equation}
where the drift $\alpha:\RR\rightarrow\RR$ is a given function.
We denote by $\cal{L}$ the infinitesimal operator of the process $(X_t)$.  
That is
\begin{equation} \label{ol} 
{\cal{L}}f :=\frac 12\partial_{xx}f -\alpha \partial_x f, \end{equation}
 We also denote by $\cal{L}^*$ the formal adjoint operator of $\cal{L}$ 
with respect to the Lebesgue measure.
In other words 
\begin{equation} \label{oal} 
{\cal{L}^*}f:=\frac 12\partial_{xx}f +\partial_x(\alpha f).  
\end{equation}

The hitting time of zero $T_0^X$ is defined as the first time that 
the process $(X_t)$ reaches zero. 
That is $T_0^X:=\inf \left\{ t\geq 0:X_t=0\right\}$.
A similar notation will be used for the hitting time of zero 
for $(B_t)$: $T_0^B$.

A probability measure $\nu$ is said to be a quasi-stationary 
distribution ($qsd$) if
\begin{equation} \label{qsd}
\forall t, \forall A \in {\cal{B}}(\RR) \ 
P_\nu (X_t\in A \mid T_0^X>t)=\nu(A).
\end{equation}
If the drift is regular, this condition is equivalent 
to the existence of $\lambda \in \RR_+$
such that 
\begin{equation} \label{qsd2}
\forall t, \forall A \in {\cal{B}}(\RR) \ 
P_\nu (X_t\in A, T_0^X>t)=e^{-\lambda t}\nu(A),
\end{equation}
which implies that $\nu=\nu_\lambda$, the probability measure concentrated
on $\RR_+$ with smooth density proportional to the unique solution
of the differential problem
\begin{equation} \label{edoqs} 
\left\{ \begin{array}{l} \cal{L}^{*}\varphi _\lambda =
-\lambda \varphi _\lambda \\ \varphi _\lambda (0)=0\quad ,\quad 
\varphi _\lambda ^{\prime }(0)=1.  \end{array} \right.  
\end{equation}
That is $\forall A \in {\cal B}(\RR_+)$ we have
$\nu_\lambda(A):=\left(\int_0^\infty\varphi_\lambda(x)dx\right)^{-1} 
\int_A\varphi_\lambda(x)dx$.

We remark that in particular from (\ref{qsd2}) the absorption
time is exponentially distributed when the initial distribution is $\nu_{\lambda}$, that is
$P_{\nu_\lambda} (T_0^X>t)=e^{-\lambda t}$.
Usually the set of values of $\lambda$ for which (\ref{edoqs}) has a
positive and integrable solution is an interval 
$(0,\underline\lambda]$ and
moreover $\underline\lambda$ coincides with the 
ground state of $\cal{L}^{*}$.
We shall prove this result in the context of a Ornstein-Uhlenbeck 
process, but this holds in many others situations 
(see \cite{Servet0}, \cite{Jaime3} and \cite{Jaime2} for example). 

This paper deals with the domain
of attraction of the $qsd$. We say that $\mu$ is in the domain
of attraction of the $qsd$ $\nu$ if
\begin{equation} \label{lqsd} \stackunder{t\rightarrow \infty }
{\lim }P_\mu (X_t\in \cdot \mid T_0^X>t)=\nu(\cdot),  
\end{equation}
where the limit is taken in the weak topology. We notice that from 
(\ref{qsd}) $\nu$ belongs
to its own domain of attraction. We work with
absolutely continuous initial distributions $\mu$, whose density $f$
is a regularly varying function (see definitions on section 3). 
Our main result is the following theorem, where we assume that
the drift is linear $\alpha=a x$ with $a>0$, that is $(X_t)$ is
an Ornstein-Uhlenbeck process.


{\bf Theorem}\hfill\break
Let $f:\RR_+\to\RR_+$ be a density function with exponent $-(1+\eta)$ 
where $\eta\in(0,1)$. 
Then for all $A\in{\cal B}(\RR_+)$
\[ \lim_{t\to\infty} P_f(X_t\in A\mid T^X_0>t)=\nu_{\lambda^*} (A)\]
where $\lambda^*:=a\eta\in(0,a)$.

\bigskip

We remark that, in the case that the density function $f$ is smooth, the condition
$\lim\limits_{u\to\infty} {uf'(u)\over f(u)}=-(1+\eta)$
is enough to ensure that $f$ has the desired exponent.

In the literature there are two works directly related with 
the problem mentioned above.  The first one was published by P. Mandl 
\cite{Mandal} who consider general drift assumptions.
He proved that the bottom of the spectrum of $\cal{L}^*$ is given by
\underline{$\lambda $}$=\sup\{\lambda\in\RR:\varphi_\lambda$ 
does not change sign$\}$. Also he proved that under certain hypothesis
on the behaviour of the Fourier transform of the initial density function $f$ 
around the point \underline{$\lambda$}, the limit in (\ref{lqsd}) 
exists and $\nu$ is 
the probability measure $\nu_{\underline \lambda}$.
The measure $\nu_{\underline{\lambda }}$ is called the minimal $qsd$.  
Because of this, we say that Mandl's result only deals with the domain 
of attraction of the minimal $qsd$. For the Ornstein-Uhlenbeck process it is
not hard to see that $\underline \lambda=a$.

The other related result is due to S.Martinez et al \cite{Jaime2}.  
They studied the domain of attraction of $qsd$  for a Brownian Motion with 
constant drift $\alpha(x)\equiv a$ with $a>0$. In this setting 
${\underline{\lambda}}=a^2/2$. They prove that the limit in 
(\ref{lqsd}) exists if $\stackunder{u\rightarrow\infty }
{\lim }-\frac{\ln f(u)}u=\beta > 0$.  
Moreover, when $\beta \in (0,a]$ $\nu$ is the probability 
measure $\nu_{\lambda^*}$ where $\lambda^*=a\beta-\frac{\beta^2}{2}$.  
When $\beta\ge a$ then $\nu=\nu_{\underline{\lambda}}$.

The present paper is organized as follows.  In section 2 we present general 
and well known facts about the Ornstein-Uhlenbeck process which we include 
for the sake of completeness.  In section 3 we present the proof of 
the main result of this paper. 


\section{Some facts about the Ornstein-Uhlenbeck process}

From now on, $(X_t)$ denotes an Ornstein-Uhlenbeck process, that is to say, 
the process that solves the stochastic differential equation (\ref{gede}) 
for a linear drift $\alpha(x):=ax$ with $a>0$ constant.  

As usual $P_x(B_t\in dy)/dy:={1\over\sqrt{2\pi t}}e^{\frac{(y-x)^2}{2t}}$ denotes 
the transition density of a Brownian Motion starting from $x$.  
We define the functions $h$ and $g$ by the formulas 
\begin{equation} \label{fhfg}
h(t):=\frac{1-e^{-2at}}{2a}\qquad,\qquad g(t):=\frac{e^{2at}-1}{2a}.  
\end{equation} 
The transition density of the process $(X_t)$ starting from $x$ can 
be computed as 
\begin{equation} 
\label{td} p(t,x,y):=P_x(X_t\in dy)/dy=P_{e^{-at}x}(B_{h(t)}\in dy)/dy
\end{equation}
On the other hand, an Ornstein-Uhlenbeck process satisfies a reflection principle.  
In words, conditioning on $T^X_0$, using the strong Markov property and the fact that 
$(X_t)$ and $(-X_t)$ have the same law under $P_0$, one obtains for $x,t>0$ 
and $A\in{\cal B}(\RR_+)$ 
\begin{equation} \label{rp1} 
P_x(X_t\in A ,T_0^X>t)=P_x(X_t\in A)-P_x(-X_t\in A).  
\end{equation}

Consequently, using formulas (\ref{td}) and (\ref{rp1}), and the reflection principle 
for a Brownian Motion, we obtain for every $x,t>0$, 
\begin{equation}
\label{rp2} P_x(T_0^X>t)=P_{e^{-at}x}(T_0^B>h(t)).
\end{equation}
We also have a formula for $q$ the transition density of the submarkovian process
$(X_t 1_{T_0^X>t})$ given by
\begin{equation}
\label{inutil}
q(t,x,y)=p(t,x,y)-p(t,x,-y)=\sqrt{2\over\pi h(t)} e^{-\{{(e^{-at}x)^2
\over2h(t)}+{y^2\over2h(t)}\}}sinh\left(e^{-at}xy\over h(t)\right),
\end{equation}
where it is assumed that $xy>0$.

We remark that formulas (\ref{fhfg}) and (\ref{rp2}) will help us 
to rewrite probabilities about the process $(X_t)$ in terms of the probabilities 
about the Brownian motion $(B_t)$.  Also, before ending the section, we notice that 
the reflection principle for the process $(X_t)$ essentially holds for any 
diffusion process $(Y_t)$ that starting from zero has the same distribution as $(-Y_t)$.  
Certainly, this is the case when a process $(Y_t)$ solves an stochastic differential 
equation, of the type (\ref{gede}), for which the drift $\alpha$ 
is an odd function and uniqueness in distribution holds.

\section{Proof of the main result}

Our main theorem relies on the concept of regularly varying functions 
(for a complete reference, see \cite{BGT}, \cite{Feller}). 
A non-negative measurable function 
$f:\RR_+\to\RR_+$ it is said to be {\bf regularly varying with exponent} 
$\bf\beta$ (briefly, $f$ has exponent $\beta$), if for all $c>0$
\[ \lim_{u\to\infty}{f(cu)\over f(u)}=c^{\beta}. \]
In order to make a discussion in the same terms of our main result, 
we deal with a function $f$ with exponent of the form $-(1+\eta)$ with $\eta\in\RR$. 
Since the function $f(u)u^{1+\eta}$ varies slowly, we have 
the following asymptotic for $\ln f$ ( see \cite{BGT}  
Proposition 1.3.6) 
\begin{equation}
\label{H1}
 \lim_{u\to\infty }\frac{\ln f(u)}{\ln u}= -(1+\eta).
\end{equation}

We start the proof of the main Theorem by proving that the set of measures \break
$\left\{P_f(X_t\in \cdot \mid T_0^X>t)\right\}_{t\ge1}$ 
is tight. The next Lemma is technically important for that purpose.

\begin{lemma}
\label{ldm}\hfill

If $f:\RR_+\to\RR_+$ is a density function with exponent $-(1+\eta)$ for some $\eta\in(0,1)$ 
then for all $\gamma\in(0,\eta)$
\begin{equation}
\label{cldm}
\lim\limits_{u\to\infty}\:u^{1-\gamma}\:{\int_u^\infty\:x^\gamma f(x)
dx\over\int_0^u xf(x)dx} = {1-\eta\over\eta-\gamma}<\infty.
\end{equation}
\end{lemma}

{\bf Proof:} From the hypothesis assumed on $f$ it follows immediately that
$x^\gamma f(x)$ is integrable near $\infty$. Therefore from \cite{Feller},
Theorem 1 on section VIII.9, we have the following limits exist
\[
\lim_{u\to\infty } \frac{u^{\gamma+1} f(u)}{\int^\infty_u x^\gamma f(x) dx}
=\eta-\gamma>0,
\]
and
\[
\lim_{u\to\infty } \frac{u^2 f(u)}{\int^u_0 x f(x) dx}
=1-\eta>0,
\]
from which the result follows.\hfill$\Box$

\begin{lemma}
\label{ldlt}\hfill

If $f:\RR_+\to\RR_+$ is a density function with 
exponent $-(1+\eta)$ for some $\eta\in(0,1)$, 
then the set of probability measures:
$\:\left\{\:P_f\left(X_t\in\cdot\:\left\vert\right.\:T^X_0>t\right)\:\right\}_{t\ge1}\:$ 
is tight.
\end{lemma}
{\bf Proof:} It is enough to check that 
$\limsup\limits_{t\to\infty}\:E_f\left(X_t^\gamma\left\vert\right.T^X_0>t\right)<\infty$ 
for some $1\ge \gamma>0$, where
\begin{equation}\label{machupichu}
E_f\left(X_t^\gamma\left\vert\right.T^X_0>t\right)
={I(0,\infty,0,\infty) \over J(0,\infty,0,\infty)},
\end{equation}
and 
\[
I(A,B,C,D):=\int_A^B\int_C^D f(x)\:y^\gamma\:q(t,x,y)\:dy\:dx,\:
J(A,B,C,D):=\int_A^B\int_C^D f(x)\:q(t,x,y)\:dy\:dx.
\]

It can be checked using (\ref{inutil}) 
that there exists a constant $c_1=c_1(a)>0$ such that for every $t\ge1$ 
and $x\in[0,e^{at}]$
\[
1\:\le\:{\int_0^\infty y^\gamma q(t,x,y)dy\over \int_0^1 
y^\gamma q(t,x,y)dy}\:\le\:c_1.
\]
From (\ref{machupichu}) and the last inequality we have that

\[
\begin{array}{lll}
 E_f\left(X_t^\gamma\left\vert\right.T^X_0>t\right) \le
{I(0,e^{at},0,\infty)\over J(0,e^{at},0,1)}
\left\{1+{I(e^{at},\infty,0,\infty) \over I(0,e^{at},0,\infty)}\right\} 
&\le& c_1{I(0,e^{at},0,1) \over J(0,e^{at},0,1)}
\left\{1+{I(e^{at},\infty,0,\infty) \over I(0,e^{at},0,1)}\right\} \\
&\le&  c_1\left\{1+{I(e^{at},\infty,0,\infty) \over I(0,e^{at},0,1)} \right\}.
\end{array}
\]

Therefore to prove the Lemma is enough to find 
$1\ge \gamma>0$ such that 

\begin{equation}
\label{nacanacalapirinaca}
\limsup\limits_{t\to\infty}{\int_{e^{at}}^\infty\int_0^\infty f(x)y^\gamma\:q(t,x,y)\:dydx
\over \int_0^{e^{at}}\int_0^1 f(x)y^\gamma\:q(t,x,y)\:dydx}\:<\:\infty.
\end{equation}

But since $0<h(1)\le h(t)\le{1\over 2a}$ for every $t\ge1$
and $sinh(z)\ge z$ for all $z\ge0$, we see from (\ref{inutil}) that

\[
\int_0^{ e^{at}}\int_0^1 f(x)y^\gamma\:q(t,x,y)\:dydx \ge
\:{4a^{3/2} \over \sqrt{\pi}e^{at}} 
\int_0^{e^{at}}xf(x)e^{-{(e^{-at}x)^2\over2h(t)}}dx
\int_0^1y^{\gamma+1}\:e^{-{y^2\over2h(t)}}dy.
\]
On the other hand $q(t,x,y)\le p(t,x,y)$. 
Thus, using (\ref{td}) and the last inequality, it follows 
the existence of a constant $c_2=c_2(a)>0$ such that

\begin{equation}
\label{tamoslistongo}
{\int_{ e^{at}}^\infty\int_0^\infty f(x)y^\gamma\:q(t,x,y)\:dydx
\over\int_0^{ e^{at}}\int_0^1 f(x)y^\gamma\:q(t,x,y)\:dydx}\le
c_2\:e^{at} \:{\int_{e^{at}}^\infty\:f(x)
\left(\int_0^\infty\:y^\gamma\:e^{-a(e^{-at}x-y)^2}dy\right)dx\over\int_0^{e^{at}}xf(x)dx}.
\end{equation}

Now, if $0< \gamma \le 1$ there exists a constant $c_3>0$ such that for $x\ge e^{at}$
\[
\int_0^\infty\:y^\gamma\:e^{-a(e^{-at}x-y)^2}dy\le c_3(e^{-at}x)^\gamma,
\]
hence from (\ref{tamoslistongo}) we can conclude that

\begin{eqnarray*}
\limsup_{t\to\infty}{\int_{ e^{at}}^\infty\int_0^\infty f(x)y^\gamma\:q(t,x,y)\:dydx
\over\int_0^{ e^{at}}\int_0^1 f(x)y^\gamma\:q(t,x,y)\:dydx} &
\le& c_2c_3\limsup_{t\to\infty} e^{at(1-\gamma)}{\int_{e^{at}}^\infty\:f(x)x^\gamma dx
\over\int_0^{e^{at}}xf(x)dx} \\
 &=& c_2c_3\limsup\limits_{u\to\infty}\:u^{1-\gamma}\:{\int_u^\infty\:x^\gamma f(x)
dx\over\int_0^u xf(x)dx}.
\end{eqnarray*}
To finish the proof we notice that Lemma \ref{ldm} ensures that the right hand side in the 
last inequality is finite for any $0<\gamma<\eta$. This proves assertion 
(\ref{nacanacalapirinaca}), and therefore the desired result follows.\hfill$\Box$

\begin{lemma}
\label{ldps}\hfill

Let $0<b<c$ and $x>0$. Then  \[1\le{P_x(T^B_0>b)\over P_x(T^B_0>c)}\le\left({c\over b}\right)^{3/2}.\]
\end{lemma}

{\bf Proof:} The first inequality is direct. On the other hand
 (see \cite{Karatzas}, page 197) after a linear substitution we get
\begin{eqnarray*}
 {P_x(T^B_0>b)\over P_x(T^B_0>c)}=
\frac{\int\limits_b^\infty {1\over\sqrt{u^3}}e^{-{x^2\over 2u}}du} 
{\int\limits_c^\infty {1\over\sqrt{u^3}}e^{-{x^2\over 2u}}du}\le
\frac{\int\limits_b^\infty {1\over\sqrt{u^3}}e^{-{x^2\over 2u}}du} 
{\int\limits_b^\infty {1\over\sqrt{(1+{c-b\over u})^3}} 
{1\over\sqrt{u^3}}e^{-{x^2\over 2u}}du}\le \left({c\over b}\right)^{3/2}.
\end{eqnarray*}

The last inequality follows since the function $u\to\left(1+{c-b\over u}\right)$ is decreasing.\hfill$\Box$

\begin{lemma}
\label{ldpj}\hfill

If $f:\RR_+\to\RR_+$ is a density function with exponent $\:-(1+\eta)$ 
where $\eta\in(0,1)$, then for every $t>0$ and $c>0$
\begin{equation}
\label{ufuf}
\lim\limits_{u\to\infty}{\int_0^{ln(u)/c}f(x)P_{c{x\over u}}(T^B_0>t)dx\over 
\int_{ln(u)/c}^\infty f(x)P_{c{x\over u}}(T^B_0>t)dx}=0.
\end{equation}
\end{lemma}

{\bf Proof:} Since the function $x\to P_x(T^B_0>t)$ is increasing on $\RR_+$
we obtain
\[ \int_0^{ln(u)/c} f(x)P_{c{x\over u}}(T^B_0>t)dx \le
P_{ln(u)\over u}(T^B_0>t) 
\int_0^{ln(u)/c} f(x)dx \le P_{ln(u)\over u}(T^B_0>t).\]

But the function $x\to P_x(T^B_0>t)$ is differentiable at $x=0$, 
hence, there exists a constant $c_1>0$, which depends on t, such that for $u>0$ 
sufficiently large
\begin{equation}
\label{pelaito}
\int_0^{ln(u)/c} f(x)P_{c{x\over u}}(T^B_0>t)dx\le c_1 {ln(u)\over u}.
\end{equation}
On the other hand, since $f$ is regularly varying, from (\ref{H1}) 
if $\kappa\in(\eta,1)$ we have for large $x$ the inequality $f(x)\ge{1\over x^{1+\kappa}}$. 
In particular, for $u>0$ large enough we have that
\[ \int_{ln(u)/c}^\infty f(x)P_{c{x\over u}}(T^B_0>t)dx\ge 
\left(c\over u\right)^\kappa \int_{ln(u)/u}^\infty {1\over x^{1+\kappa}}P_x(T^B_0>t)dx.\]
Now, since the function $u\to {ln(u)\over u}$ is asymptotically decreasing to $0$, 
there exists a constant $c_2>0$, which also depends on $t$, such that for big $u>0$
\begin{equation}
\label{pelaote}
\int_{ln(u)/c}^\infty f(x)P_{c{x\over u}}(T^B_0>t)dx\ge {c_2\over u^\kappa}.
\end{equation}
From (\ref{pelaito}) and (\ref{pelaote}), and the fact that $\kappa<1$, 
it follows the result.\hfill$\Box$

\begin{lemma}
\label{ldpf/pf}\hfill

If $f:\RR_+\to\RR_+$ is a density function with exponent $\:-(1+\eta)$ 
where $\eta\in(0,1)$, then for every $s\ge0$

\begin{equation}
\label{manu}
\lim\limits_{t\to\infty}{P_f(T^X_0>t+s)\over P_f(T^X_0>t)}=e^{-(a\eta)s}.
\end{equation}
\end{lemma}

{\bf Proof:} Let $s>0$ and $h$ be the function defined on (\ref{fhfg}). 
Then, using (\ref{rp2}), it follows for every $t>0$ that

\[
{P_f(T^X_0>t+s)\over P_f(T^X_0>t)}={\int_0^\infty f(x) 
P_{e^{-a(t+s)}x}(T^B_0>h(t+s))dx\over\int_0^\infty f(x) P_{e^{-at}x}(T^B_0>h(t))dx}.
\]

Notice that $h(t)<h(t+s)<{1\over 2a}$ and 
$\lim\limits_{t\to\infty}h(t)={1\over 2a}$. 
Therefore, using the bound obtained in Lemma \ref{ldps} and setting $u=e^{at}$, 
it can be easily seeing that (\ref{manu}) is equivalent to

\begin{equation}
\label{ifeas/ifeas}
\lim_{u\to\infty} {\int_0^\infty f(x)
\label{if2a/if2a}
    P_{e^{-as}{x\over u}}(T^B_0>{1\over 2a})dx\over\int_0^\infty f(x)
    P_{x\over u}(T^B_0>{1\over 2a})dx}=e^{-(a\eta)s}.
\end{equation}

Writing

\[
{\int_0^\infty f(x) P_{e^{-as}{x\over u}}(T^B_0>{1\over 2a})dx
\over\int_0^\infty f(x) P_{x\over u}(T^B_0>{1\over 2a})dx}=
{\left(\int\limits_0^{e^{as}ln(u)}+
\int\limits_{e^{as}ln(u)}^\infty\right) f(x) P_{e^{-as}{x\over u}}(T^B_0>{1\over 2a})dx
\over \left(\int\limits_0^{ln(u)}\quad+\quad
\int\limits_{ln(u)}^\infty\right)\quad f(x) P_{x\over u}(T^B_0>{1\over 2a})dx},
\]

from lemma \ref{ldpj} and using the substitution $y=e^{-as}x$ in the numerator, 
to check (\ref{ifeas/ifeas}) it is sufficient to prove  that
\[
\lim_{u\to\infty}{\int_{ln(u)}^\infty f(e^{as}x) P_{x\over
u}(T^B_0>{1\over 2a})dx\over \int_{ln(u)}^\infty f(x) P_{x\over
u}(T^B_0>{1\over 2a})dx }=e^{-a(1+\eta)s}.
\]
But, $f$ is regularly varying therefore $\lim\limits_{x\to\infty}{f(e^{as}x)\over f(x)}
=e^{-a(1+\eta)s}$ for all $x>0$. From this fact, it follows the result.\hfill$\Box$

\begin{lemma}\hfill
\label{ldfqs}

\begin{description}
\item[(a)] $\left\{\lambda\:\vert\:\varphi_\lambda\hbox{ does not change sign }\right\}=(-\infty,a]$;
\item[(b)] For every $\lambda\in(0,a]$, $\int_0^\infty\varphi_\lambda(x)dx<\infty$.
\end{description}
\end{lemma}

{\bf Proof:} Let $\varphi_\lambda$ be the solution of (\ref{edoqs}). Then  
$\psi_\lambda(u)=e^{a u^2}\varphi_\lambda(u)$ is the unique solution of the 
equation
\[
\label{eqspfpsi}
\left\{
\begin{array}{l}
\mathrel{{\cal L}\psi_\lambda}\:=\:-\lambda\:\psi_\lambda \\
\psi_\lambda(0)=0 \quad,\quad\psi_\lambda'(0)=1.
\end{array}
\right.
\]
Hence, to prove \ref{ldfqs}(a), we just need to concentrate our attention on $\psi_\lambda$. 
Note that $\psi_\lambda$ is an analytic function. 
Setting $\psi_\lambda(u)=\sum\limits_{k\ge0}b_{k} u^{k}$ one sees that 
$b_0=0,\ b_1=1$ and one obtains a recursion
for $b_{k+2}$ in terms of $b_k$ which shows that  $\forall k\ge 0 \ b_{2k}=0$ and for $k \ge 1$
\begin{eqnarray}
\label{csm1}
b_{2k+1}={a^k\over (2k+1) k!}\prod\limits_{i=0}^{k-1}
\left( 1-{\lambda a^{-1}\over 2i+1} \right).
\end{eqnarray}

In particular $\psi_\lambda$ cannot change sign if $\lambda\le a$. Notice also that 
$\psi_a(u)=u$.  

On the other hand, if $\lambda\in(a,3a)$ then from (\ref{csm1}) we get that $b_{2k+1}<0$ 
for every $k\ge 1$. Thus, $\lim\limits_{u\to\infty}\psi_\lambda(u)=-\infty$. 
But $\psi_\lambda(0)=0$ and $\psi_\lambda'(0)>0$, hence, from the last limit, 
we see that there exists $x_0>0$ such that $\psi_\lambda>0$ on $(0,x_0)$ and 
$\psi_\lambda(x_0)=0$. 

Let $\lambda>a$. We prove then that $\psi_\lambda$ has to change its sign. 
Letting $\kappa\in(a,min\{3a\:,\lambda\})$ we just proved the 
existence of some $x_0>0$ such that $\psi_\kappa>0$ on $(0,x_0)$ 
and $\psi_\kappa(x_0)=0$. But, simultaneously

\[
\left\{
\begin{array}{l}
\left(e^{-au^2}\psi_\lambda'(u)\right)'+2\lambda e^{-au^2}\psi_\lambda(u)=0 \\
\left(e^{-au^2}\psi_\kappa'(u)\right)'+2\kappa e^{-au^2}\psi_\kappa(u)=0.
\end{array}
\right.
\]

Since $\lambda>\kappa$, by the Sturn-Liouville's theorem 
(see \cite{Sotomayor}, page 104), there exists $y_0\in(0,x_0)$ 
such that $\psi_\lambda(y_0)=0$. 
Hence, to prove that $\psi_\lambda$ changes sign, it is 
sufficient to check that $\psi_\lambda'(y_0)\ne0$. If this were not 
the case it can be easily checked out that $\psi_\lambda^{(k)}(y_0)=0$ for all $k\ge0$. 
Thus, $\psi_\lambda\equiv0$ on $[0,\infty)$. This can not occur because 
$\psi_\lambda'(0)\ne0$. Consequently, for $\lambda>a$, $\psi_\lambda'(y_0)\ne0$
 which implies that $\psi_\lambda$ changes its sign. This proves \ref{ldfqs}(a).

Now, we prove part \ref{ldfqs}(b). For $\lambda\in(0,a]$, we have 
\begin{eqnarray}
\int_0^\infty\varphi_\lambda(u)du &=& {1\over 2a}\sum_{k\ge0} \:{k!\over a^k}\:b_{2k+1} 
={1\over 2a}\left\{1+\sum_{k\ge1} {1\over2k+1}\prod_{i=0}^{k-1}
\left( 1-{\lambda a^{-1}\over 2i+1} \right)\right\} \nonumber \\
&\le& {1\over 2a}\left\{1+\sum_{k\ge1} \:{1\over 2k+1}\:exp
\left(-\lambda a^{-1}\sum_{i=0}^{k-1}{1\over 2i+1}\right)\right\} \\
&\le& {1\over 2a}\left\{1+\sum_{k\ge1} 
\:{1\over 2k+1}\:exp\left(-{\lambda \over 2a}\sum_{i=1}^{k}{1\over i}\right)\right\}.
\label{tdn}
\end{eqnarray}
where we have used the fact $1-z\le e^{-z}$, 
for $z\in[0,1]$. Let $d:=\sup\limits_{k>0}\left\vert 
ln(k)-\sum\limits_{i=1}^k{1\over i}
\right\vert$. Since $d<\infty$, (\ref{tdn}) yields
\[
\int_0^\infty\varphi_\lambda(u)du \le
{1\over 2a}\left\{1+{e^{\lambda d\over 2a} \over 2} \sum_{k\ge1}
{1\over k^{1+{\lambda \over 2a}}}\right\}.
\]

The previous inequality shows that  $\int_0^\infty\varphi_\lambda(u)du<\infty$, 
which proves \ref{ldfqs}(b).\hfill$\Box$

\bigskip
{\bf Proof of the Theorem}

Let $t_n'\to\infty$. From lemma \ref{ldlt}, we know that there exits a subsequence 
$t_n\to\infty$ and a probability measure $\mu$, defined on ${\cal B}(\RR_+)$, such that
\begin{equation}
\label{wlim}
\lim\limits_{n\to\infty}P_f(X_{t_n}\in \cdot \mid T_0^X>t_n)=\mu.
\end{equation}

Now, the function $P_x(T_0^X>s)$ is continuous and bounded on $x\in\RR_+$.
The strong Markov property allows us to show that for every $n$

\[
\frac{P_f(T_0^X>t_n+s)}{P_f(T_0^X>t_n)}= E_f\left(P_{X_{t_n}}(T_0^X>s)\mid T_0^X>t_n\right).
\]

Therefore, taking limit as $n\to\infty$ and using lemma \ref{ldpf/pf}, 
we deduce that for every $s>0$
\begin{equation}
\label{mami}
\int_0^\infty P_x(T_0^X>s)\mu (dx)=e^{-(a\eta)s}=e^{-\lambda^*s},
\end{equation}
where $\lambda^*:=a\eta\in(0,a)$. 
On the other hand, by lemma \ref{ldfqs}(b), we see that 
$\varphi_{\lambda^*}\ge0$ is 
integrable on $[0,\infty)$. 
This allows us to define for each $s>0$ and $y>0$ the function
\[
\Lambda(s,y):=\int_0^\infty q(s,x,y) \varphi_{\lambda^*}(x)dx.
\]

A simple computation shows
\[
\frac{\partial}{\partial s}q(s,x,y)={\cal L}_x q(s,x,y),
\]
which is nothing but the classical Kolmogorov's equation. 
Now, from the given definitions we have

\begin{eqnarray*}
\frac\partial{\partial s}\Lambda (s,y)&=&
\int_0^\infty {\cal L}_x q(s,x,y) \varphi_{\lambda^*} (x)dx 
=\int_0^\infty q(s,x,y) ({\cal L}_x^{*}\varphi_{\lambda^*}) (x)dx \\
&=&-\lambda^*\int_0^\infty q(s,x,y) \varphi_{\lambda^*} (x)dx 
\label{edos1}= -\lambda^* \Lambda(s,y).
\end{eqnarray*}
Therefore $\Lambda(s,y)=\Lambda(0,y) e^{-\lambda^*s}=
\varphi_{\lambda^*}(y) e^{-\lambda^*s}$.
Hence, integrating this equality over $y$ leads to
\begin{equation}
\label{idpl}
\int_0^\infty\!\!\int_0^\infty q(s,x,y) dy 
\:\varphi_{\lambda^*}(x)dx=\int_0^\infty P_x(T_0^X>s) 
\;\varphi_{\lambda^*}(x) (dx)
=c e^{-\lambda^* s},
\end{equation}
where $c=\int_0^\infty\varphi_{\lambda^*}(x)dx \in (0,\infty)$.

Now, define on ${\cal B}(\RR_+)$ the finite signed-measure 
$\rho(A):=\mu(A)-{1\over c}\int_A\varphi_{\lambda^*}(x)dx$.
It follows that $\rho([0,\infty))=0$. 
From (\ref{mami}) and (\ref{idpl}), and then (\ref{rp2}), 
we see that for every $s>0$

\begin{equation}
\label{completa}
\int_0^\infty P_{e^{-as}x}(T^B_0>h(s))\rho(dx)=0.
\end{equation}

Using the well known distribution for the running maxima
of a Brownian motion one has $P_x(T^B_0>t)={2\over\sqrt{2\pi}}
\int_0^{x\over\sqrt{t}} e^{-{y^2\over2}}dy$.
This identity, a simple substitution and integrating by parts allow 
us to conclude from (\ref{completa}) that for every $s>0$

\[
\int_0^\infty \rho[0,x] \; {1\over \sqrt{2\pi g(s)}}e^{-{x^2\over 2g(s)}}dx=0, 
\]

where $g$ is defined as in (\ref{fhfg}). Finally, because of the set of 
density functions 
$\left\{{1\over \sqrt{2\pi\theta}}e^{-{x^2\over \theta}}\right\}_{\theta>0}$ 
is a complete family and the range of $g$ is $[0,\infty)$,
we deduce that 
$\rho[0,x]=0$, for every $x>0$. Therefore, for all $A\in{\cal B}(\RR_+)$,
 $\mu(A)={1\over c}\int_A\varphi_{\lambda^*} (x)dx=\nu_{\lambda^*}(A)$. 
Notice that the limiting 
measure does not depend on the initially chosen sequence, 
from which the result follows.\hfill$\Box$

\bigskip
{\bf Acknowledgments.}
We are grateful to an anonymous referee for her/his valuable
suggestions and comments, and for pointing out references \cite{BGT}
and \cite{Feller}.
The authors are indebted to 
C\'atedra Presidencial fellowship, FONDAP program on Applied
Mathematics and project FONDEF.
\bigskip

\end{document}